\RequirePackage{lineno}
\documentclass[12pt,a4paper]{article}

\usepackage[left=2cm,right=2cm, top=2cm,bottom=2cm,bindingoffset=0cm]{geometry}

\usepackage{color}

\usepackage{subfig}
\usepackage{flushend}
\usepackage{indentfirst}
\usepackage{graphics}
\usepackage{amsmath}
\usepackage{graphicx}
\usepackage{epstopdf}
\usepackage{float}
\usepackage{amssymb}
\usepackage[font=footnotesize,labelfont=bf]{caption}
\usepackage[labelsep=period]{caption}

\usepackage[font=footnotesize,labelfont=bf]{caption}
\usepackage[labelsep=period]{caption}

\graphicspath{{./figure/}}

\newtheorem{thm}{Theorem}[section]
\newtheorem{lem}{Lemma}[section]
\newtheorem{ozn}{Definition}
\newtheorem{nas}{Corollary}[section]

\begin{document}

\title{{Filtering of periodically correlated processes}}

\date{}

\maketitle

\noindent Applied Statistics. Actuarial and Financial Mathematics.  No. 2, 149--158,  2012\\
ISSN 1810-3022

\vspace{20pt}

\author{\textbf{Iryna Dubovets'ka}, \textbf{Mykhailo Moklyachuk}$^*$    \\\\
\emph{ Department of Probability Theory, Statistics and Actuarial
Mathematics, \\
Taras Shevchenko National University of Kyiv, Kyiv 01601, Ukraine}\\
Corresponding email: moklyachuk@gmail.com}\\\\\\

\noindent \textbf{\large{Abstract}} \hspace{2pt}
The problem of optimal linear estimation of a linear functional depending on the unknown values of periodically correlated stochastic process from observations of the process with additive noise is considered. 
Formulas for calculating the  mean square error and the spectral characteristic of the optimal linear estimate of the functional are proposed in the case where spectral densities are exactly known.
Formulas that determine the least favorable spectral densities and the minimax (robust) spectral characteristics are proposed for a given class of admissible spectral densities.
\\

\noindent\textbf{{Keywords}} \hspace{2pt}
periodically correlated process, robust estimate, mean square error, least favorable spectral density, minimax spectral characteristic.
\\

\noindent\textbf{ AMS 2010 subject classifications.} Primary: 60G25, 60G35 Secondary: 62M20, 93E10

\noindent\hrulefill

\section{{Introduction}}

The study of mean-square continuous periodically correlated processes was initiated in the article by E. G. Gladyshev [1], where properties of the correlation function and representations of periodically correlated processes were analyzed. The connection between periodically correlated and stationary processes was studied by A. Makagon [2], [3].
Due to the correspondence relations between the processes, the problem of estimation of periodically correlated processes is reduced to the corresponding problem for vector stationary sequences.
The main results regarding the representations of periodically correlated sequences through simpler random sequences are presented in the work of L. Hurd and A. Miami [4].

Methods for studying the problems of estimation of unknown values of stationary processes (extrapolation, interpolation, and filtration problems)
were developed in the works by A.N. Kolmogorov [5], N. Wiener [6], A.M. Yaglom [7], [8].
The developed methods are based on the assumption that the exact values of the spectral densities of the processes are known.
In the case where complete information on the spectral densities is impossible, while a set of admissible spectral densities is given,
the minimax method of solution the estimation problems is used. That is, an estimate is found that minimizes the value of the error simultaneously for all densities from a given class.
U. Grenander [9] first applied the minimax approach to the problem of extrapolation of stationary processes.
M. P. Moklyachuk [10-15], M. P. Moklyachuk, and O. Yu. Masyutka [16] investigated the problems of extrapolation, interpolation,
and filtering for stationary processes and sequences.
Minimax problems of optimal estimation of linear functionals from periodically correlated sequences and processes were studied in the works of
I. I. Dubovets'ka, O. Yu. Masyutka, and M. P. Moklyachuk [17-21].

In this article the problem of mean-square optimal linear estimation of the functional $A\zeta =\int_{0}^{\infty }{a(t)\zeta (-t)\,dt}$ from unknown values of the mean-square continuous periodically correlated process $\zeta (t)$ based on the results of observations of the process $\zeta (t)+\theta (t)$ at points $t\le 0$, where $\theta (t)$ is a periodically correlated process uncorrelated with $\zeta (t)$, is investigated. Formulas for calculating the spectral characteristic and the mean-square error of the optimal estimate of the functional $A\zeta $ are presented.
For given classes of admissible spectral densities, the least favorable spectral densities and the minimax spectral characteristic of the optimal linear estimate of the functional $A\zeta $  are determined.

\section{Periodically correlated processes and corresponding vector stationary sequences}

\begin{ozn} \label{def1} [1] A mean-square continuous stochastic process $\zeta :\mathbb{R}\to H={{L}_{2}}(\Omega ,F,P)$, $E\zeta (t)=0$,
is called periodically correlated (PC) with period $T$ if its correlation function $K(t+u,u)=E\zeta (t+u)\overline{\zeta (u)}$ for all $t,u\in \mathbb{R}$ and some fixed $T>0$ satisfies the condition
$K(t+u,u)=K(t+u+T,u+T).$
\end{ozn}

Let $\{\zeta (t),\,t\in \mathbb{R}\}$ and $\{\theta (t),\,t\in \mathbb{R}\}$ be uncorrelated PC stochastic processes.
Let us construct two sequences of stochastic functions
\begin{equation}\label{eq01}
     \left\{ {{\zeta }_{j}}(u)=\zeta (u+jT),\,u\in [0,T),\,j\in \mathbb{Z} \right\},
\end{equation}
\begin{equation}\label{eq02}
     \left\{ {{\theta }_{j}}(u)=\theta (u+jT),\,u\in [0,T),\,j\in \mathbb{Z} \right\}.
\end{equation}

Each of the sequences (1), (2) forms  ${{L}_{2}}([0,T);H)$-valued stationary sequence $\{{{\zeta }_{j}},\,j\in \mathbb{Z}\}$ and $\{{{\theta }_{j}},\,j\in \mathbb{Z}\}$,
respectively, with correlation functions
\begin{multline*}
{{B}_{\zeta }}(l,j)={{\left\langle {{\zeta }_{l}},{{\zeta }_{j}} \right\rangle }_{H}}=\int_{o}^{T}{E\,\zeta (u+lT)}\overline{\zeta (u+jT)}\,du=\\
=\int_{o}^{T}{{{K}_{\zeta }}(u+(l-j)T,u)}\,du={{B}_{\zeta }}(l-j),
     \end{multline*}
     \begin{multline*}
    {{B}_{\theta }}(l,j)={{\left\langle {{\theta }_{l}},{{\theta }_{j}} \right\rangle }_{H}}=
    \int_{o}^{T}{E\,\theta (u+lT)}\overline{\theta (u+jT)}\,du=\\
    =\int_{o}^{T}{{{K}_{\theta }}(u+(l-j)T,u)}\,du={{B}_{\theta }}(l-j),
    \end{multline*}
where ${{K}_{\zeta }}(t,s)=E\zeta (t)\overline{\zeta (s)},$ ${{K}_{\theta }}(t,s)=E\theta (t)\overline{\theta (s)}$ are the correlation functions of PC processes $\zeta (t)$ and $\theta (t)$.

If we define an orthonormal basis in ${{L}_{2}}([0,T);\mathbb{R})$
\[\left\{{{\tilde{e}}_{k}}=\frac{1}{\sqrt{T}}{{e}^{2\pi i\{{{(1)}^{-k}}\left[ {}^{k}/{}_{2} \right]\}u/T}},\,k=1,2,...\right\},\quad\left\langle {{{\tilde{e}}}_{k}},{{{\tilde{e}}}_{j}} \right\rangle ={{\delta }_{jk}},\]
then stationary sequences $\{{{\zeta }_{j}},\,j\in \mathbb{Z}\}$, $\{{{\theta }_{j}},\,j\in \mathbb{Z}\}$
 can be represented in the form
\begin{equation}\label{eq03}
     {{\zeta }_{j}}=\sum\nolimits_{k=1}^{\infty }{{{\zeta }_{kj}}{{{\tilde{e}}}_{k}}}, \,\,  {{\zeta }_{kj}}=\left\langle {{\zeta }_{j}},{{{\tilde{e}}}_{k}} \right\rangle =\frac{1}{\sqrt{T}}\int_{0}^{T}{{{\zeta }_{j}}(v){{e}^{-2\pi i\{{{(1)}^{-k}}\left[ {}^{k}/{}_{2} \right]\}v/T}}}\,dv,
\end{equation}
\begin{equation}\label{eq04}
{{\theta }_{j}}=\sum\nolimits_{k=1}^{\infty }{{{\theta }_{kj}}{{{\tilde{e}}}_{k}}},\,\,\,{{\theta }_{kj}}=\left\langle {{\theta }_{j}},{{{\tilde{e}}}_{k}} \right\rangle =\frac{1}{\sqrt{T}}\int_{0}^{T}{{{\theta }_{j}}(v){{e}^{-2\pi i\{{{(1)}^{-k}}\left[ {}^{k}/{}_{2} \right]\}v/T}}}\,dv.
\end{equation}

The components ${{\zeta }_{kj}}$ and ${{\theta }_{kj}}$ of the stationary sequences $\{{{\zeta }_{j}},\,j\in \mathbb{Z}\}$ and $\{{{\theta }_{j}},\,j\in \mathbb{Z}\}$ satisfy the conditions [10, 22]

     \[E{{\zeta }_{kj}}=0,\,\,   ||{{\zeta }_{j}}||_{H}^{2}=\sum\nolimits_{k=1}^{\infty }{E|{{\zeta }_{kj}}{{|}^{2}}}={{P}_{\zeta }}<\infty,\,\,\,
     E{{\zeta }_{kl}}\overline{{{\zeta }_{nj}}}=\left\langle {{R}_{\zeta }}(l-j){{e}_{k}},{{e}_{n}} \right\rangle ,\]
     \[E{{\theta }_{kj}}=0,\,\,\,||{{\theta }_{j}}||_{H}^{2}=\sum\nolimits_{k=1}^{\infty }{E|{{\theta }_{kj}}{{|}^{2}}}={{P}_{\theta }}<\infty ,\,\,\,E{{\theta }_{kl}}\overline{{{\theta }_{nj}}}=\left\langle {{R}_{\theta }}(l-j){{e}_{k}},{{e}_{n}} \right\rangle .\]
where $\{{{e}_{k}},\,\,k=1,2,...\}$ is a basis in the space ${{\ell }_{2}}$.
The correlation functions ${{R}_{\zeta }}(j)$ and ${{R}_{\theta }}(j)$ of stationary sequences $\{{{\zeta }_{j}},\,j\in \mathbb{Z}\}$ and $\{{{\theta }_{j}},\,j\in \mathbb{Z}\}$ are operator functions in ${{\ell }_{2}}$.
Correlation operators ${{R}_{\zeta }}(0)={{R}_{\zeta }}$, ${{R}_{\theta }}(0)={{R}_{\theta }}$ are kernel operators and
\[\sum\nolimits_{k=1}^{\infty }{\left\langle {{R}_{\zeta }}{{e}_{k}},{{e}_{k}} \right\rangle }=\,||{{\zeta }_{j}}||_{H}^{2}={{P}_{\zeta }},
\quad
\sum\nolimits_{k=1}^{\infty }{\left\langle {{R}_{\theta }}{{e}_{k}},{{e}_{k}} \right\rangle }=\,||{{\theta }_{j}}||_{H}^{2}={{P}_{\theta }}.\]

Stationary sequences $\{{{\zeta }_{j}},\,j\in \mathbb{Z}\}$ and $\{{{\theta }_{j}},\,j\in \mathbb{Z}\}$ have spectral densities
$f(\lambda )=\{{{f}_{kn}}(\lambda )\}_{k,n=1}^{\infty },$ $g(\lambda )=\{{{g}_{kn}}(\lambda )\}_{k,n=1}^{\infty },$
which are positive operator-valued functions in ${{\ell }_{2}}$ of the variable $\lambda \in [-\pi ,\pi )$,
if their correlation functions ${{R}_{\zeta }}(j)$ and ${{R}_{\theta }}(j)$ can be represented in the form
\[\left\langle {{R}_{\zeta }}(j){{e}_{k}},{{e}_{n}} \right\rangle =\frac{1}{2\pi }\int_{-\pi }^{\pi }{{{e}^{ij\lambda }}\left\langle f(\lambda ){{e}_{k}},{{e}_{n}} \right\rangle }\,d\lambda
\]
\[\left\langle {{R}_{\theta }}(j){{e}_{k}},{{e}_{n}} \right\rangle =\frac{1}{2\pi }\int_{-\pi }^{\pi }{{{e}^{ij\lambda }}\left\langle g(\lambda ){{e}_{k}},{{e}_{n}} \right\rangle }\,d\lambda ,\,\,k,n=1,2,....\]

For almost all $\lambda \in [-\pi ,\pi )$ the spectral densities $f(\lambda ),$ $g(\lambda )$ are kernel operators with integrable kernel norms:
     \[\sum\nolimits_{k=1}^{\infty }{\frac{1}{2\pi }\int_{-\pi }^{\pi }{\left\langle f(\lambda ){{e}_{k}},{{e}_{k}} \right\rangle }\,d\lambda =\sum\nolimits_{k=1}^{\infty }{\,\left\langle {{R}_{\zeta }}{{e}_{k}},{{e}_{k}} \right\rangle }}=||{{\zeta }_{j}}||_{H}^{2}={{P}_{\zeta }},\]
     \[\sum\nolimits_{k=1}^{\infty }{\frac{1}{2\pi }\int_{-\pi }^{\pi }{\left\langle g(\lambda ){{e}_{k}},{{e}_{k}} \right\rangle }\,d\lambda =\sum\nolimits_{k=1}^{\infty }{\,\left\langle {{R}_{\theta }}{{e}_{k}},{{e}_{k}} \right\rangle }}=||{{\theta }_{j}}||_{H}^{2}={{P}_{\theta }}.\]

\section{Classic optimal linear filtering method}

We study the problem of mean-square optimal linear estimation of the functional
$$A\zeta =\int_{0}^{\infty }{a(t)\zeta (-t)\,dt}$$
from unknown values of the mean-square continuous PC process $\zeta (t)$ based on observations of the process $\zeta (t)+\theta (t)$ at points $t\le 0$, where $\zeta (t)$ is an uncorrelated PC process with $\theta (t)$. The function $a(t),\,t\in {{\mathbb{R}}_{+}},$ satisfies the condition $\int_{0}^{\infty }{|a(t)|\,dt}<\infty .$

Let us write the functional $A\zeta $ in the following form
     $$A\zeta =\int_{0}^{\infty }{a(t)\zeta (-t)\,dt=}\sum\nolimits_{j=0}^{\infty }{\int_{0}^{T}{{{a}_{j}}(u){{\zeta }_{-j}}(-u)}\,du,}$$
     $${{a}_{j}}(u)=a(u+jT),\quad {{\zeta }_{-j}}(-u)=\zeta (-u-jT),\quad u\in [0,T).$$

Taking in account the decomposition (3) of the stationary sequence $\{{{\zeta }_{j}},\,j\in \mathbb{Z}\}$, the functional $A\zeta $ can be represented as follows
$$A\zeta =\sum\nolimits_{j=0}^{\infty }{\int_{0}^{T}{{{a}_{j}}(u){{\zeta }_{-j}}(-u)}\,du=\sum\nolimits_{j=0}^{\infty }{\sum\nolimits_{k=1}^{\infty }{{{a}_{kj}}{{\zeta }_{k,-j}}}=\sum\nolimits_{j=0}^{\infty }{\vec{a}_{j}^{\top }{{{\vec{\zeta }}}_{-j}}}},}$$
where
$${{\vec{\zeta }}_{-j}}={{({{\zeta }_{k,-j}},\,k=1,2,...)}^{\top }},$$
$$\,\,{{\vec{a}}_{j}}={{({{a}_{kj}},\,k=1,2,...)}^{\top }}={{({{a}_{1j}},{{a}_{3j,}}{{a}_{2j}},...,{{a}_{2k+1,j}},{{a}_{2k,j}},...)}^{\top }},$$
$${{a}_{kj}}=\left\langle {{a}_{j}},{{{\tilde{e}}}_{k}} \right\rangle =\frac{1}{\sqrt{T}}\int_{0}^{T}{{{a}_{j}}(v){{e}^{-2\pi i\{{{(1)}^{-k}}\left[ {}^{k}/{}_{2} \right]\}v/T}}}\,dv.$$
Assume that the coefficients $\{{{\vec{a}}_{j}},\,j=0,1,...\}$ satisfy the conditions
\begin{equation}\label{eq05}
\sum\nolimits_{j=0}^{\infty }{||{{{\vec{a}}}_{j}}||}<\infty,\quad {{\sum\nolimits_{j=0}^{\infty }{(j+1)||{{{\vec{a}}}_{j}}||}}^{2}}<\infty,\quad ||{{\vec{a}}_{j}}|{{|}^{2}}=\sum\nolimits_{k=1}^{\infty }{|{{a}_{kj}}{{|}^{2}}}.
\end{equation}

\begin{ozn} \label{def2} [22]
Denote by ${{H}_{\zeta }}(n)$  a closed linear subspace of the Hilbert space $H$ generated by random variables $\{{{\zeta }_{kj}},\,k\ge 1\,,\,j\le n\}.$ The sequence $\{{{\zeta }_{j}},\,j\in \mathbb{Z}\}$ is called regular if $\bigcap\nolimits_{n}{{{H}_{\zeta }}(n)=}\emptyset.$
  If $\bigcap\nolimits_{n}{{{H}_{\zeta }}(n)=}H$, then the sequence $\{{{\zeta }_{j}},\,j\in \mathbb{Z}\}$ is called singular.
 \end{ozn}

Since the unknown values of the components of a singular stationary sequence are estimated without error, we can consider the problem of optimal linear estimation only for regular stationary sequences.

A regular stationary sequence $\{{{\zeta }_{j}}+{{\theta }_{j}},\,j\in \mathbb{Z}\}$ admits a canonical representation of the moving average of its components [10, 22]
 \begin{equation}\label{eq06}
       {{\zeta }_{kj}}+{{\theta }_{kj}}=\sum\nolimits_{u=-\infty }^{j}{\sum\nolimits_{m=1}^{M}{{{d}_{km}}(j-u){{\varepsilon }_{m}}(u)}},
  \end{equation}
where ${{\varepsilon }_{m}}(u),\,m=1,...,M,\,u\in \mathbb{Z}$ are mutually orthogonal sequences in $H$ with orthogonal values:
$E\,{{\varepsilon }_{m}}(u)\overline{{{\varepsilon }_{p}}(v)}={{\delta }_{mp}}{{\delta }_{uv}};$ $M$ is the multiplicity of the stationary sequence
$\{{{\zeta }_{j}},\,j\in \mathbb{Z}\};$ sequences ${{d}_{km}}(u),\,\,k=1,2,...,\,\,m=1,...,\,M,\,u=0,1,...,$
are such that
$$
\sum_{u=0}^{\infty }\sum_{k=1}^{\infty }\sum_{m=1}^{M}|{d}_{km}|^{2}<\infty.
$$
The spectral density of such a stationary sequence $\{{{\zeta }_{j}}+{{\theta }_{j}},\,j\in \mathbb{Z}\}$ admits the canonical factorization
 \begin{equation}\label{eq07}
f(\lambda )+g(\lambda )=P(\lambda )\,{{P}^{*}}(\lambda ),\quad P(\lambda )=\sum\nolimits_{u=0}^{\infty }{d(u){{e}^{-iu\lambda }}},
\end{equation}
where the matrix $d(u)=\{{{d}_{km}}(u)\}_{k=\overline{1,\infty }}^{m=\overline{1,M}}$ is determined by the coefficients of the canonical representation (6).

The spectral density $f(\lambda )$ admits canonical factorization if
 \begin{equation}\label{eq08}
f(\lambda )=\varphi (\lambda )\,{{\varphi }^{*}}(\lambda ),\quad \varphi (\lambda )=\sum\nolimits_{u=0}^{\infty }{\varphi (u){{e}^{-iu\lambda }}},
\end{equation}
where $\varphi (\lambda )=\{{{\varphi }_{km}}(\lambda )\}_{k=\overline{1,\infty }}^{m=\overline{1,M}}.$

The spectral density $g(\lambda )$ admits canonical factorization if
 \begin{equation}\label{eq09}
g(\lambda )=\psi (\lambda )\,{{\psi }^{*}}(\lambda ),\quad \psi (\lambda )=\sum\nolimits_{u=0}^{\infty }{\psi (u){{e}^{-iu\lambda }}},
\end{equation}
where $\psi (\lambda )=\{{{\psi }_{km}}(\lambda )\}_{k=\overline{1,\infty }}^{m=\overline{1,M}}.$ To factorize (7) the density $f(\lambda )+g(\lambda )$, it is sufficient to factorize one of the densities (8) or (9).

Let us denote by ${{L}_{2}}(f)$ the Hilbert space of vector functions $b(\lambda )=\{{{b}_{k}}(\lambda )\}_{k=1}^{\infty }$, which are integrated in measure with the density $f(\lambda )$:
$$\int_{-\pi }^{\pi }{{{b}^{\top }}(\lambda )f(\lambda )}\overline{b(\lambda )}\,d\lambda =\int_{-\pi }^{\pi }{\sum\nolimits_{k,n=1}^{\infty }{{{b}_{k}}(\lambda ){{f}_{kn}}(\lambda )\overline{{{b}_{n}}(\lambda )}}}\,d\lambda <\infty .$$

Let us denote by $L_{2}^{-}(f)$ the subspace in ${{L}_{2}}(f)$ generated by the functions of the form ${{e}^{ij\lambda }}{{\delta }_{k}},$ $j\le 0,$${{\delta }_{k}}=\{{{\delta }_{kn}}\}_{n=1}^{\infty },$ $k=1,2,...,$ where ${{\delta }_{kn}}$ is the Kronecker symbol: ${{\delta }_{kk}}=1$, ${{\delta }_{kn}}=0$ for $k\ne n$.

 The linear estimate $\hat{A}\zeta $ of the functional $A\zeta $ based on the observations of the sequence $\{{{\zeta }_{j}}+{{\theta }_{j}}\}$ at points $j\le 0$ is determined by the spectral characteristic $h({{e}^{i\lambda }})\in L_{2}^{-}(f+g)$ and has the form
 \begin{equation}\label{eq10}
\hat{A}\zeta =\int_{-\pi }^{\pi }{{{h}^{\top }}({{e}^{i\lambda }})}({{Z}^{\zeta +\theta }}(\,d\lambda ))=\int_{-\pi }^{\pi }{\sum\nolimits_{k=1}^{\infty }{{{h}_{k}}({{e}^{i\lambda }})(Z_{k}^{\zeta +\theta }(\,d\lambda )}}),
\end{equation}
where ${{Z}^{\zeta +\theta }}(\Delta )=\{Z_{k}^{\zeta +\theta }(\Delta )\}_{k=1}^{\infty }$ is an orthogonal random measure of the sequence $\{{{\zeta }_{j}}+{{\theta }_{j}}\}$.

The mean square error of the linear estimate $\hat{A}\zeta $ with the spectral characteristic $h({{e}^{i\lambda }})=\sum\nolimits_{j=0}^{\infty }{{{{\vec{h}}}_{j}}{{e}^{-ij\lambda }}}$ can be calculated by the formula
 $$\Delta (h;f,g)=E|A\zeta -\hat{A}\zeta {{|}^{2}}=$$
 $$=\frac{1}{2\pi }\int_{-\pi }^{\pi }{{{[A({{e}^{i\lambda }})]}^{\top }}g(\lambda )\overline{A({{e}^{i\lambda }})}}\,d\lambda +\frac{1}{2\pi }\int_{-\pi }^{\pi }{{{[A({{e}^{i\lambda }})-h({{e}^{i\lambda }})]}^{\top }}(f(\lambda )+g(\lambda ))\overline{[A({{e}^{i\lambda }})-h({{e}^{i\lambda }})]}}\,d\lambda -$$
 $$-\frac{1}{2\pi }\int_{-\pi }^{\pi }{{{[A({{e}^{i\lambda }})-h({{e}^{i\lambda }})]}^{\top }}g(\lambda )\overline{A({{e}^{i\lambda }})}}\,d\lambda -\frac{1}{2\pi }\int_{-\pi }^{\pi }{{{[A({{e}^{i\lambda }})]}^{\top }}g(\lambda )\overline{[A({{e}^{i\lambda }})-h({{e}^{i\lambda }})]}}\,d\lambda =$$
     $$=\,||\Psi a|{{|}^{2}}+||D(a-h)|{{|}^{2}}-\left\langle \Psi (a-h),\Psi a \right\rangle -\left\langle \Psi a,\Psi (a-h) \right\rangle ,$$
where the action of the operators $\Psi $ and $D$ is given as follows
$$A({{e}^{i\lambda }})=\sum\nolimits_{j=0}^{\infty }{{{{\vec{a}}}_{j}}{{e}^{-ij\lambda }}},\quad ||\Psi a|{{|}^{2}}=\sum\nolimits_{q=0}^{\infty }{||{{(\Psi a)}_{q}}|{{|}^{2}},}\quad {{(\Psi a)}_{q}}=\sum\nolimits_{l=0}^{q}{{{\psi }^{\top }}(q-l){{{\vec{a}}}_{l}}},$$
$$||D(a-h)|{{|}^{2}}=\sum\nolimits_{q=0}^{\infty }{||{{(D(a-h))}_{q}}|{{|}^{2}},}\quad {{(D(a-h))}_{q}}=\sum\nolimits_{l=0}^{q}{{{d}^{\top }}(q-l)\,({{{\vec{a}}}_{l}}-{{{\vec{h}}}_{l}})},$$
 $$
 \left\langle \Psi (a-h),\Psi a \right\rangle =\overline{\left\langle \Psi a,\Psi (a-h) \right\rangle }=\sum\nolimits_{q=0}^{\infty }{\left\langle {{(\Psi (a-h))}_{q}},{{(\Psi a)}_{q}} \right\rangle }.
 $$
The spectral characteristic $h(f,g)$ of the optimal linear estimate $\hat{A}\zeta $ with given densities $f(\lambda ),\,g(\lambda )$ minimizes the mean square error
 \begin{equation}\label{eq11}
      \Delta (f,g)=\Delta (h(f,g);f,g)=\underset{h\in L_{2}^{-}(f,g)}{\mathop{\min }}\,\Delta (h;f,g)=\underset{\hat{A}\zeta }{\mathop{\min }}\,E\,|A\zeta -\hat{A}\zeta {{|}^{2}}.
  \end{equation}

Suppose that the densities $f(\lambda )+g(\lambda )$ and $g(\lambda )$ admit factorizations (7) and (9). Then the spectral characteristic $h(f,g)$, which is the solution of problem (11), and the mean square error $\Delta (f,g)$ of the optimal estimate $\hat{A}\zeta $ are calculated by the formulas
 \begin{equation}\label{eq12}
 h(f,g)=A({{e}^{i\lambda }})-{{b}^{\top }}(\lambda ){{S}_{g}}({{e}^{i\lambda }}),
   \end{equation}
 \begin{equation}\label{eq13}
 \Delta (f,g)=||\Psi a|{{|}^{2}}-||{{B}^{*}}{{\Psi }^{*}}\Psi a|{{|}^{2}},
   \end{equation}
where
    $$b(\lambda )=\sum\nolimits_{u=0}^{\infty }{b(u){{e}^{-iu\lambda }}},\quad b(\lambda )d(\lambda )={{I}_{M}},$$
 $$
    {{S}_{g}}({{e}^{i\lambda }})=\sum\nolimits_{l=0}^{\infty }{{{({{S}_{g}})}_{l}}{{e}^{-il\lambda }}},\quad {{({{S}_{g}})}_{l}}={{({{B}^{*}}{{\Psi }^{*}}\Psi a)}_{l}}=\sum\nolimits_{j=0}^{\infty }{\overline{b(j)}}{{({{\Psi }^{*}}\Psi a)}_{l+j}},$$
    \begin{equation}\label{eq14}
     {{({{\Psi }^{*}}\Psi a)}_{j}}=\sum\nolimits_{u=0}^{\infty }{\overline{\psi (u)}{{(\Psi a)}_{u+j}}},\quad ||{{B}^{*}}{{\Psi }^{*}}\Psi a|{{|}^{2}}=\sum\nolimits_{l=0}^{\infty }{||({{B}^{*}}{{\Psi }^{*}}\Psi a}{{)}_{l}}|{{|}^{2}}.
      \end{equation}

 If the densities $f(\lambda )+g(\lambda )$ and $f(\lambda )$ allow factorizations (7) and (8), then the spectral characteristic $h(f,g)$ and the mean square error $\Delta (f,g)$ of the optimal estimate $\hat{A}\zeta $ are calculated by the formulas
    \begin{equation}\label{eq15}
    h(f,g)={{b}^{\top }}(\lambda ){{S}_{f}}({{e}^{i\lambda }}),
     \end{equation}
       \begin{equation}\label{eq16}
     \Delta (f,g)=||\Phi a|{{|}^{2}}-||{{B}^{*}}{{\Phi }^{*}}\Phi a|{{|}^{2}},
      \end{equation}
where
     $${{S}_{f}}({{e}^{i\lambda }})=\sum\nolimits_{l=0}^{\infty }{{{({{S}_{f}})}_{l}}{{e}^{-il\lambda }}},\quad {{({{S}_{f}})}_{l}}={{({{B}^{*}}{{\Phi }^{*}}\Phi a)}_{l}}=\sum\nolimits_{j=0}^{\infty }{\overline{b(j)}}{{({{\Phi }^{*}}\Phi a)}_{l+j}},$$
$${{({{\Phi }^{*}}\Phi a)}_{j}}=\sum\nolimits_{u=0}^{\infty }{\overline{\varphi (u)}{{(\Phi a)}_{u+j}}},\quad {{(\Phi a)}_{q}}=\sum\nolimits_{l=0}^{q}{{{\varphi }^{\top }}(q-l){{{\vec{a}}}_{l}}},$$
   \begin{equation}\label{eq17}
||\Phi a|{{|}^{2}}=\sum\nolimits_{q=0}^{\infty }{||{{(\Phi a)}_{q}}|{{|}^{2}},}\quad ||{{B}^{*}}{{\Phi }^{*}}\Phi a|{{|}^{2}}=\sum\nolimits_{l=0}^{\infty }{||({{B}^{*}}{{\Phi }^{*}}\Phi a}{{)}_{l}}|{{|}^{2}}.
 \end{equation}

 Therefore, the following theorem is valid.

\begin{thm} \label{thm1}
Let $\{\zeta (t),\,t\in \mathbb{R}\}$ and $\{\theta (t),\,t\in \mathbb{R}\}$ be uncorrelated PC random processes such that the stationary sequences $\{{{\zeta }_{j}},\,j\in \mathbb{Z}\}$ and $\{{{\theta }_{j}},\,j\in \mathbb{Z}\}$, constructed according to relations (1), (2), respectively, have spectral densities $f(\lambda )$ and $g(\lambda )$, which admit canonical factorizations (7), (8) or (7), (9). Let the coefficients $\{{{\vec{a}}_{j}},j=0,1,...\},$ that define the functional $A\zeta $ satisfy conditions (5). Then the spectral characteristic $h(f,g)$ and the mean square error $\Delta (f,g)$ of the estimate of the functional $A\zeta $ from observations of the process $\zeta (t)+\theta (t)$ at points $t\le 0$ are calculated by formulas (15), (16) or (12), (13), respectively.
The linear optimal estimate of the functional $A\zeta $ is determined by formula (10).
\end{thm}

\begin{nas} \label{corr1}
Let $\{\zeta (t),\,t\in \mathbb{R}\}$ and $\{\theta (t),\,t\in \mathbb{R}\}$ be uncorrelated PC random processes. Let one of the stationary sequences $\{{{\zeta }_{j}},\,j\in \mathbb{Z}\}$ or $\{{{\theta }_{j}},\,j\in \mathbb{Z}\}$ constructed according to relations (1), (2), respectively, be a vector sequence of white noise with a coordinate dispersion ${{\sigma }^{2}}$. Then the spectral characteristic $h(f,g)$ of the optimal linear estimate of the functional $A\zeta $ is calculated by formula (15) or (12). The mean square error of the estimate equals
     $$\Delta (f,g)={{\sigma }^{2}}||\vec{a}|{{|}^{2}}-{{\sigma }^{4}}||{{B}^{*}}a|{{|}^{2}},$$
where
 $$\vec{a}=({{\vec{a}}_{j}},\,j=0,1,...),\, ||\vec{a}|{{|}^{2}}=\sum\nolimits_{j=0}^{\infty }{||{{{\vec{a}}}_{j}}}|{{|}^{2}},$$ $$ ||{{B}^{*}}a|{{|}^{2}}=\sum\nolimits_{l=0}^{\infty }{||{{({{B}^{*}}a)}_{l}}|{{|}^{2}}},\, {{({{B}^{*}}a)}_{l}}=\sum\nolimits_{j=0}^{\infty }{\overline{b(j)\,}{{{\vec{a}}}_{l+j}}}.$$
\end{nas}

\begin{nas} \label{corr2}  Under the conditions of Corollary \ref{corr1}, the  mean square error of the optimal linear estimate of the value $\vec{a}_{N}^{\top }{{\vec{\zeta }}_{-N}}$ based on  observations of the process $\zeta (t)+\theta (t)$ at points $t\le 0$ is calculated by the formula
$$\Delta (f,g)={{\sigma }^{2}}||{{\vec{a}}_{N}}|{{|}^{2}}-{{\sigma }^{4}}\sum\nolimits_{q=0}^{N}{||\overline{b(q)}{{{\vec{a}}}_{N}}|{{|}^{2}}}.$$
 \end{nas}

 \section{{Minimax (robust) filtering}}
 To use formulas (12), (13), (15), (16) to calculate the spectral characteristic and the mean square error of the optimal estimate of the functional $A\zeta $, it is necessary to know the spectral densities $f(\lambda )$ and $g(\lambda )$ of the stationary sequences $\{{{\zeta }_{j}},\,j\in \mathbb{Z}\}$ and $\{{{\theta }_{j}},\,j\in \mathbb{Z}\}$, which are constructed according to relations (1), (2), respectively. In the case where the spectral densities are not exactly known, while a set $D={{D}_{f}}\times {{D}_{g}}$ of admissible spectral densities is given, the minimax approach is resonable to use for solution  the problems of estimating functionals from unknown process values. We find an estimate that gives the smallest error simultaneously for all spectral densities from a given class $D$.

\begin{ozn} \label{def3}
 For a given set of pairs of spectral densities $D={{D}_{f}}\times {{D}_{g}}$ the spectral densities ${{f}^{0}}(\lambda )\in {{D}_{f}},$ ${{g}^{0}}(\lambda )\in {{D}_{g}}$ are called the least favorable in $D$ for the optimal estimation of the functional $A\zeta $ if
$$\Delta ({{f}^{0}},{{g}^{0}})=\Delta (h({{f}^{0}},{{g}^{0}});{{f}^{0}},{{g}^{0}})=\underset{(f,g)\in D}{\mathop{\max }}\,\Delta (h(f,g);f,g).$$
 \end{ozn}

\begin{ozn} \label{def4}
For a given set of pairs of spectral densities $D={{D}_{f}}\times {{D}_{g}}$ the spectral characteristic ${{h}^{0}}(\lambda )$ of the optimal estimate of the functional $A\zeta $ is called minimax (robust) if
$${{h}^{0}}(\lambda )\in {{H}_{D}}=\bigcap\nolimits_{(f,g)\in D}{L_{2}^{-}(f+g)},\quad \underset{h\in {{H}_{D}}}{\mathop{\min }}\,\underset{(f,g)\in D}{\mathop{\max }}\,\Delta (h;f,g)=\underset{(f,g)\in D}{\mathop{\max }}\,\Delta ({{h}^{0}};f,g).$$
 \end{ozn}

Taking into account these definitions and the above relations (7)-(13), we can verify that the following lemmas hold true.

\begin{lem} \label{lem1}
The spectral densities ${{f}^{0}}(\lambda )\in {{D}_{f}}$ and ${{g}^{0}}(\lambda )\in {{D}_{g}}$, which admit canonical factorizations (7)-(9),
will be the least favorable in the class $D$ for the optimal estimate of $A\zeta $ if the coefficients of the factorizations (7)-(9) determine the solutions of the conditional extremum problem
 \begin{equation}\label{eq18}
\Delta (f,g)=||\Phi a|{{|}^{2}}-||{{B}^{*}}{{\Phi }^{*}}\Phi a|{{|}^{2}}\to \,\sup ,
 \end{equation}
$$f(\lambda )=\left( \sum\nolimits_{u=0}^{\infty }{\varphi (u){{e}^{-iu\lambda }}} \right){{\left( \sum\nolimits_{u=0}^{\infty }{\varphi (u){{e}^{-iu\lambda }}} \right)}^{*}}\in {{D}_{f}},$$
$$g(\lambda )=\left( \sum\nolimits_{u=0}^{\infty }{d(u){{e}^{-iu\lambda }}} \right){{\left( \sum\nolimits_{u=0}^{\infty }{d(u){{e}^{-iu\lambda }}} \right)}^{*}}-\left( \sum\nolimits_{u=0}^{\infty }{\varphi (u){{e}^{-iu\lambda }}} \right){{\left( \sum\nolimits_{u=0}^{\infty }{\varphi (u){{e}^{-iu\lambda }}} \right)}^{*}}\in {{D}_{g}},$$
or conditional extremum problem
\begin{equation}\label{eq19}
\Delta (f,g)=||\Psi a|{{|}^{2}}-||{{B}^{*}}{{\Psi }^{*}}\Psi a|{{|}^{2}}\to \,\sup ,
 \end{equation}
$$g(\lambda )=\left( \sum\nolimits_{u=0}^{\infty }{\psi (u){{e}^{-iu\lambda }}} \right){{\left( \sum\nolimits_{u=0}^{\infty }{\psi (u){{e}^{-iu\lambda }}} \right)}^{*}}\in {{D}_{g}},$$
$$f(\lambda )=\left( \sum\nolimits_{u=0}^{\infty }{d(u){{e}^{-iu\lambda }}} \right){{\left( \sum\nolimits_{u=0}^{\infty }{d(u){{e}^{-iu\lambda }}} \right)}^{*}}-\left( \sum\nolimits_{u=0}^{\infty }{\psi (u){{e}^{-iu\lambda }}} \right){{\left( \sum\nolimits_{u=0}^{\infty }{\psi (u){{e}^{-iu\lambda }}} \right)}^{*}}\in {{D}_{f}}.$$
 \end{lem}

\begin{lem} \label{lem2}
 Let the spectral density $f(\lambda )$ be known and admits the canonical factorization (8). Then the spectral density ${{g}^{0}}(\lambda )$ admits the canonical factorizations (7), (9) and is the least favorable for the optimal estimate of the functional $A\zeta $ if
$$f(\lambda )+{{g}^{0}}(\lambda )=\left( \sum\nolimits_{u=0}^{\infty }{{{d}^{0}}(u){{e}^{-iu\lambda }}} \right){{\left( \sum\nolimits_{u=0}^{\infty }{{{d}^{0}}(u){{e}^{-iu\lambda }}} \right)}^{*}},$$
where the coefficients $\{{{d}^{0}}(u),\,u=0,1,...\}$ are determined by the solutions $\{{{d}^{0}}(u),\,u=0,1,...\}$ of the conditional extremum problem
\begin{equation}\label{eq20}
||{{B}^{*}}{{\Phi }^{*}}\Phi a|{{|}^{2}}\to \,\inf,\quad g(\lambda )=\left( \sum\nolimits_{u=0}^{\infty }{d(u){{e}^{-iu\lambda }}} \right){{\left( \sum\nolimits_{u=0}^{\infty }{d(u){{e}^{-iu\lambda }}} \right)}^{*}}-f(\lambda )\in {{D}_{g}}.
 \end{equation}
 \end{lem}

\begin{lem} \label{lem3}
 Let the spectral density $g(\lambda )$ be known and admits the canonical factorization (9). Then the spectral density ${{f}^{0}}(\lambda )$ admits the canonical factorizations (7), (8) and is the least favorable for the optimal estimate of the functional $A\zeta $ if
$${{f}^{0}}(\lambda )+g(\lambda )=\left( \sum\nolimits_{u=0}^{\infty }{{{d}^{0}}(u){{e}^{-iu\lambda }}} \right){{\left( \sum\nolimits_{u=0}^{\infty }{{{d}^{0}}(u){{e}^{-iu\lambda }}} \right)}^{*}},$$
where the coefficients $\{{{d}^{0}}(u),\,u=0,1,...\}$ are determined by the solutions $\{{{d}^{0}}(u),\,u=0,1,...\}$ of the conditional extremum problem
\begin{equation}\label{eq21}
||{{B}^{*}}{{\Psi }^{*}}\Psi a|{{|}^{2}}\to \,\inf,\quad f(\lambda )=\left( \sum\nolimits_{u=0}^{\infty }{d(u){{e}^{-iu\lambda }}} \right){{\left( \sum\nolimits_{u=0}^{\infty }{d(u){{e}^{-iu\lambda }}} \right)}^{*}}-g(\lambda )\in {{D}_{f}}.
 \end{equation}
 \end{lem}

The least favorable spectral densities ${{f}^{0}}(\lambda )\in {{D}_{f}}$, ${{g}^{0}}(\lambda )\in {{D}_{g}}$ and the minimax spectral characteristic ${{h}^{0}}=h({{f}^{0}},{{g}^{0}})$ form a saddle point of the function $\Delta (h;f,g)$ on the set ${{H}_{D}}\times D$. Saddle point inequalities
 $$\Delta ({{h}^{0}};f,g)\le \Delta ({{h}^{0}};{{f}^{0}},{{g}^{0}})\le \Delta (h;{{f}^{0}},{{g}^{0}}),\quad \forall h\in {{H}_{D}},\quad \forall f\in {{D}_{f}},\quad \forall g\in {{D}_{g}},$$
 are satisfied if ${{h}^{0}}=h({{f}^{0}},{{g}^{0}})$, $h({{f}^{0}},{{g}^{0}})\in {{H}_{D}}$ and $({{f}^{0}},{{g}^{0}})$ is a solution to the conditional extremum problem
\begin{equation}\label{eq22}
\Delta (h({{f}^{0}},{{g}^{0}});f,g)\to \,\sup,\quad (f,g)\in D,
 \end{equation}
where the functional
$$\Delta (h({{f}^{0}},{{g}^{0}});f,g)=$$
$$=\frac{1}{2\pi }\int_{-\pi }^{\pi }{{{[S_{g}^{0}({{e}^{i\lambda }})]}^{\top }}{{b}^{0}}(\lambda )f(\lambda ){{({{b}^{0}}(\lambda ))}^{*}}\overline{S_{g}^{0}({{e}^{i\lambda }})}}\,d\lambda +\frac{1}{2\pi }\int_{-\pi }^{\pi }{{{[S_{f}^{0}({{e}^{i\lambda }})]}^{\top }}{{b}^{0}}(\lambda )g(\lambda ){{({{b}^{0}}(\lambda ))}^{*}}\overline{S_{f}^{0}({{e}^{i\lambda }})}}\,d\lambda $$
linearly depends on the unknown densities $(f,g)$ from the set of admissible densities $D$, the functions $S_{f}^{0}({{e}^{i\lambda }}),$ $S_{g}^{0}({{e}^{i\lambda }})$ are calculated by formulas (17), (14) provided that $f(\lambda )={{f}^{0}}(\lambda ),$ $g(\lambda )={{g}^{0}}(\lambda ).$

\section{{Least favorable spectral densities in the class ${{D}_{0,0}}$}}

Let us consider the problem of minimax estimation of the functional $A\zeta $ from the PC process $\{\zeta (t),\,t\in \mathbb{R}\}$ for the set of spectral densities of vector stationary sequences $\{{{\zeta }_{j}},\,j\in \mathbb{Z}\}$ and $\{{{\theta }_{j}},\,j\in \mathbb{Z}\}$, respectively, which are constructed according to relations (1), (2):
$${{D}_{0,0}}=\left\{ (f(\lambda ),g(\lambda )|\frac{1}{2\pi }\int_{-\pi }^{\pi }{{{f}_{kk}}(\lambda )\,d\lambda }={{p}_{k}},\,\frac{1}{2\pi }\int_{-\pi }^{\pi }{{{g}_{kk}}(\lambda )\,d\lambda }={{q}_{k}},\,k=1,2,... \right\}.$$
The set ${{D}_{0,0}}$ characterizes restrictions on the second moment of the processes $\{\zeta (t),\,t\in \mathbb{R}\}$ and $\theta (t),\,t\in \mathbb{R}\}$.

Using the method of indefinite Lagrange multipliers, we find that the solution $({{f}^{0}},{{g}^{0}})$ of  the conditional extremum problem (22) satisfies the following relations:
$${{\left( {{b}^{0}}(\lambda ) \right)}^{\top }}S_{g}^{0}({{e}^{i\lambda }}){{\left( S_{g}^{0}({{e}^{i\lambda }})\, \right)}^{*}}\overline{{{b}^{0}}(\lambda )}=\left\{ \alpha _{k}^{2}{{\delta }_{kn}} \right\}_{k.n=1}^{\infty },$$
$${{\left( {{b}^{0}}(\lambda ) \right)}^{\top }}S_{f}^{0}({{e}^{i\lambda }}){{\left( S_{f}^{0}({{e}^{i\lambda }})\, \right)}^{*}}\overline{{{b}^{0}}(\lambda )}=\left\{ \beta _{k}^{2}{{\delta }_{kn}} \right\}_{k.n=1}^{\infty },$$
where $\alpha _{k}^{2},\,\,\beta _{k}^{2},\,\,k=1,2,...$ are the undetermined Lagrange multipliers. The last equations can be transformed as follows
\begin{multline}\label{eq23}
\left( \sum\nolimits_{l=0}^{\infty }{{{(S_{g}^{0})}_{l}}{{e}^{-il\lambda }}} \right){{\left( \sum\nolimits_{l=0}^{\infty }{{{(S_{g}^{0})}_{l}}{{e}^{-il\lambda }}} \right)}^{*}}=
\\
=\left( \sum\nolimits_{u=0}^{\infty }{{{d}^{0}}(u){{e}^{-iu\lambda }}} \right)\left\{ \alpha _{k}^{2}{{\delta }_{kn}} \right\}_{k.n=1}^{\infty }{{\left( \sum\nolimits_{u=0}^{\infty }{{{d}^{0}}(u){{e}^{-iu\lambda }}} \right)}^{*}},
\end{multline}
\begin{multline}\label{eq24}
\left( \sum\nolimits_{l=0}^{\infty }{{{(S_{f}^{0})}_{l}}{{e}^{-il\lambda }}} \right){{\left( \sum\nolimits_{l=0}^{\infty }{{{(S_{f}^{0})}_{l}}{{e}^{-il\lambda }}} \right)}^{*}}=
\\
=\left( \sum\nolimits_{u=0}^{\infty }{{{d}^{0}}(u){{e}^{-iu\lambda }}} \right)\left\{ \beta _{k}^{2}{{\delta }_{kn}} \right\}_{k.n=1}^{\infty }{{\left( \sum\nolimits_{u=0}^{\infty }{{{d}^{0}}(u){{e}^{-iu\lambda }}} \right)}^{*}}.
\end{multline}

The unknown Lagrange multipliers $\alpha _{k}^{2},\,\,\beta _{k}^{2},\,\,k=1,2,...$, the coefficients $\{{{b}^{0}}(u),\,u=0,1,...\}$ are determined from the canonical factorization equations (7)-(9) of the spectral densities $f(\lambda )+g(\lambda )$, $f(\lambda ),$ $g(\lambda )$ and the constraints imposed on the densities by the class ${{D}_{0,0}}$.
If one of the spectral densities is known, then one of the relations (23) or (24) is used to calculate the least favorable spectral densities of the given class ${{D}_{0,0}}$.

The following theorem is verified.

\begin{thm}\label{thm2}
  The least favorable spectral densities ${{f}^{0}}(\lambda ),\,{{g}^{0}}(\lambda )$ in the class ${{D}_{0,0}}$ for the optimal estimate of the functional $A\zeta $ are determined from equations (23) and (24), factorizations (7)-(9), from  the conditional extremum problems (18) or (19) and from restrictions of the class ${{D}_{0,0}}$.
  The minimax spectral characteristic $h({{f}^{0}},{{g}^{0}})$ of the estimate $\hat{A}\zeta $ is calculated by formula (15) or (12).
  The mean square error $\Delta ({{f}^{0}},{{g}^{0}})$ is calculated by formula (16) or (13).
\end{thm}

\begin{nas}\label{corr3} If the spectral density matrix $f(\lambda )\,$ (or $g(\lambda )$) is known and admits canonical factorization (8) (respectively (9)), then the least favorable spectral density $\,{{g}^{0}}(\lambda )$ (${{f}^{0}}(\lambda )$) is determined by relations (7) - (9), (20), (24) ((7)-(9), (21), (23)) and restrictions of the class ${{D}_{0,0}}$.
The minimax spectral characteristic $h({{f}^{0}},{{g}^{0}})$ of the estimate $\hat{A}\zeta $ is calculated by formula (15) or (12).
The mean square error $\Delta ({{f}^{0}},{{g}^{0}})$ is calculated by formula (16) or (13).
\end{nas}

\section{{Conclusions}}

 In this paper, formulas are derived for calculating the mean square error and spectral characteristic for the problem of optimal estimation of the functional
 $A\zeta =\int_{0}^{\infty }{a(t)\zeta (-t)\,dt}$ from unknown values of the mean square continuous periodically correlated process $\zeta (t)$
 based on  observations of the process $\zeta (t)+\theta (t)$ at points $t\le 0$, where $\theta (t)$ is a periodically correlated process uncorrelated with $\zeta (t)$.
 The problem is studied in the case of spectral certainty, i.e. where the spectral densities are known, and in the case of spectral uncertainty, i.e. where the spectral densities are unknown while a class of admissible spectral densities is given.
 The least favorable spectral densities and minimax (robust) spectral characteristics of optimal estimates of the functional $A\zeta $ are determined for a certain class of admissible spectral densities ${{D}_{0,0}}$.

\noindent\hrulefill

\end{document}